\documentclass[11pt,reqno]{amsart}
\usepackage{mathrsfs,eucal}
\usepackage{hyperref,amsbsy}
\usepackage{graphicx}

\newtheorem{theorem}{Theorem}[section]
\newtheorem{lemma}[theorem]{Lemma}

\newtheorem{corollary}[theorem]{Corollary}

\theoremstyle{definition}

\theoremstyle{remark}

\numberwithin{equation}{section}

\newcommand{\lra}{\longrightarrow}

\newcommand{\smap}[3]{#1 : #2 \lra #3\xspace}
\newcommand{\set}[1]{\{#1\}\xspace}

\renewcommand{\l}{\ell}

\DeclareSymbolFont{AMSb}{U}{msb}{m}{n}
\DeclareMathSymbol{\F}{\mathbin}{AMSb}{"46}
\DeclareMathSymbol{\N}{\mathbin}{AMSb}{"4E}
\DeclareMathSymbol{\Z}{\mathbin}{AMSb}{"5A}
\DeclareMathSymbol{\R}{\mathbin}{AMSb}{"52}
\DeclareMathSymbol{\C}{\mathbin}{AMSb}{"43}

\usepackage{xspace}

\begin{document}

\title[Fibonacci Identities and Graph Colorings]{Fibonacci Identities and Graph Colorings}% \\[1ex]{\small{\it Working Notes: \today}}}

\author{Christopher J. Hillar}%
\address{Department of Mathematics, Texas A\&M University, College Station, TX 77843, USA.}
\email{chillar@math.tamu.edu}
\thanks{The first author is supported under a National Science Foundation Postdoctoral Fellowship.}

%   Information for second author 
 \author{Troels Windfeldt}
\address{Department of Mathematical Sciences, University of Copenhagen, Denmark.}
\email{windfeldt@math.ku.dk}

\subjclass{05C15, 05A19, 05C38}
\keywords{Graph colorings, Fibonacci identities, chromatic polynomial, linear recurrences}

%---------------------------------------------------------
\begin{abstract}
We generalize both the Fibonacci and Lucas numbers to the context of graph colorings, and prove some identities involving these numbers. As a corollary we obtain new proofs of some known identities involving Fibonacci numbers such as 
\[F_{r+s+t} = F_{r+1}F_{s+1}F_{t+1} + F_r F_s F_t - F_{r-1}F_{s-1}F_{t-1}.\]  
\end{abstract} 
\maketitle
%---------------------------------------------------------

% General info

%\cite{latexbook}

%---------------------------------------------------------------------------
% Begin SECTION
%---------------------------------------------------------------------------
\section{Introduction}
In graph theory, it is natural to study vertex colorings, and more specifically, those
colorings in which adjacent vertices have different colors.  
In this case, the number of such colorings of a graph $G$ is encoded by 
the chromatic polynomial of $G$.  This object can be computed using the method of ``deletion and contraction",
which involves the recursive combination of chromatic polynomials for smaller graphs.
The purpose of this note is to show how the Fibonacci and Lucas numbers 
(and other integer recurrences) arise naturally in this context, 
and in particular, how identities among these numbers
can be generated from the different choices for decomposing a graph into smaller pieces.

We first introduce some notation.
Let $G$ be a undirected graph (possibly containing loops and multiple edges) with vertices 
$V=\set{1, \ldots, n} $ and edges $E$.
Given nonnegative integers $k$ and $\l$, a \emph{$(k,\l)$-coloring} of
$G$ is a map \[\smap{\varphi}{V}{\{c_1, \ldots,c_{k+\l}\}},\]
in which $\{c_1, \ldots, c_{k+\l}\}$ is a fixed set of $k+\l$ ``colors".
The map $\varphi$ is called \emph{proper} if whenever $i$ is adjacent to $j$ and $\varphi(i),\varphi(j) \in \set{c_1, \ldots, c_k}$, we have $\varphi(i) \neq \varphi(j)$.  Otherwise, we
say that the map $\varphi$ is \emph{improper}.
In somewhat looser terminology, one can think of $\{c_{k+1},\ldots,c_{k+\l}\}$ as coloring
``wildcards".

Let $\chi_G(x,y)$ be a function such that $\chi_G(k,\l)$ is the number of proper $(k,\l)$-colorings of $G$.
This object was introduced by the authors of \cite{DPT03} and 
can be given as a polynomial in $x$ and $y$ (see Lemma \ref{edgerecursion}).  
It simultaneously generalizes the chromatic, independence, and matching polynomials of $G$.
For instance, $\chi_G(x,0)$ is the usual chromatic polynomial while 
$\chi_G(x,1)$ is the independence polynomial for $G$ (see \cite{DPT03} for more
details).

We next state a simple rule that enables one to calculate the polynomial $\chi_G(x,y)$ recursively.
In what follows, $G\backslash e$ denotes the graph obtained by removing the edge $e$ from
$G$, and for a subgraph $H$ of $G$, the graph $G\backslash H$ is gotten from $G$ by 
removing $H$ and all the edges of $G$ that are adjacent to vertices of $H$.  Additionally,
the \textit{contraction} of an edge $e$ in $G$ is the graph $G/e$ obtained by removing $e$
and identifying as equal the two vertices sharing this edge.

\begin{lemma}\label{edgerecursion}
Let $e$ be an edge in $G$, and let $v$ be the vertex to which $e$ contracts in
$G/e$. Then,
\begin{equation}\label{edgerecursioneq}
\chi_{G}(x,y)=\chi_{G\backslash e}(x,y) - \chi_{G/e}(x,y) + y\cdot \chi_{(G/e)\backslash v}(x,y).
\end{equation}
\end{lemma}
\begin{proof}
The number of proper $(k,\l)$-colorings of $G\backslash e$ which 
have distinct colors for the vertices sharing edge $e$ is given by
$\chi_{G\backslash e}(k,\l) - \chi_{G/e}(k,\l)$; these colorings are also proper for $G$.
The remaining proper $ (k,\l) $-colorings of $ G $ are precisely those for which the vertices sharing edge $ e $ have the same color. This color must be one of the wildcards $\{c_{k+1},\ldots,c_{k+\l}\}$, and so the number of remaining proper $ (k,\l) $-colorings of $ G $ is counted by \mbox{$\l \cdot \chi_{(G/e)\backslash v}(k,\l)$}.
\end{proof}

With such a recurrence, we need to specify initial conditions.  When $G$ simply consists of one vertex and has no edges, we have
$\chi_{G}(x,y) = x+y$, and when $G$ is the empty graph, we set $\chi_{G}(x,y) = 1$ (consider
$G$ with one edge joining two vertices in (\ref{edgerecursioneq})).  Moreover, $\chi$ is
multiplicative on disconnected components.  This allows us to compute $\chi_G$ for
any graph recursively.

In the special case when $k = 1$, there is also a way to calculate $\chi_{G}(1,y)$ by
removing vertices from $G$.  Define the \textit{link} of a vertex $v$ to 
be the subgraph link$(v)$ of $G$ consisting of $v$, the edges touching $v$, and 
the vertices sharing one of these edges with $v$.  Also if $u$ and $v$ are joined by an edge $e$,  
we define $\textup{link}(e)$ to be $\textup{link}(u) \cup \textup{link}(v)$ in $G$, and also
we set $\deg(e)$ to be $\deg(u) + \deg(v) - 2$.  We then have the following rules.

\begin{lemma}\label{vertexrecursion}
Let $v$ be any vertex of $G$, and let $e$ be any edge.  Then,
\begin{align}
\label{vertexrecur1}
\chi_{G}(1,y) = \ & y \cdot \chi_{G \backslash v}(1,y) + y^{\deg(v)} \cdot \chi_{G \backslash  \text{\rm link}(v)}(1,y),\\ \label{vertexrecur2}
\chi_{G}(1,y) = \ & \phantom{ \, \,  \, y \cdot } \chi_{G \backslash e}(1,y) - y^{\deg(e)} \cdot \chi_{G \backslash  \text{\rm link}(e)}(1,y).
\end{align}

\end{lemma}
\begin{proof}
The number of proper $(1,\l)$-colorings of $G$ with vertex $v$ colored with a wildcard
is $\l \cdot \chi_{G \backslash v}(1,\l)$.  Moreover, in any proper coloring of $G$ with $v$ colored $c_1$, 
each vertex among the $\deg(v)$ ones adjacent to $v$ can only be one of the $\l$ wildcards.
This explains the first equality in the lemma.  

Let $v$ be the vertex to which $e$ contracts in $G/e$.  From
equation (\ref{vertexrecur1}), we have
\[  \chi_{G/e}(1,y)= y \cdot \chi_{(G/e)\backslash v}(1,y) +  y^{\deg(v)} \cdot \chi_{(G/e) \backslash \text{\rm link}(v)}(1,y).\]
Subtracting this equation from (\ref{edgerecursioneq}) with $x=1$,
and noting that $\deg(e) = \deg(v)$ and 
$G \backslash \textup{link}(e) = (G/e) \backslash \text{\rm link}(v)$,
we arrive at the second equality in the lemma.
\end{proof}

Let $P_n$ be the path graph on $n$ vertices and let $C_n$ be the cycle graph, also
on $n$ vertices ($C_1$ is a vertex with a loop attached while $ C_2 $ is two vertices joined by two edges).  
Fixing nonnegative integers $k$ and $\l$ not both zero, we define the 
following sequences of numbers ($n \geq 1$):
\begin{equation}
\begin{split}
a_{n} = \ & \chi_{P_n} (k,\l), \\
b_{n} = \ & \chi_{C_n} (k,\l). \\
\end{split}
\end{equation}
As we shall see, these numbers are natural generalizations of both the Fibonacci and 
Lucas numbers to the context of graph colorings.  The following lemma uses 
graph decomposition to give simple recurrences for these sequences.

\begin{lemma}\label{recurlem}
The sequences $a_{n}$ and $b_n$ satisfy the following linear recurrences with initial conditions:
\begin{gather}
\begin{flalign}
a_1 & = k+\l, & a_2 & = (k+\l)^2-k, & a_{n} & = (k+\l-1)a_{n-1} + \l a_{n-2}; \\
b_1 & = \l, & b_2 & = (k+\l)^2-k, & b_3 & = a_3 - b_2 +\l a_1,
\end{flalign}\\
\begin{flalign}
b_{n} & = (k+\l-2)b_{n-1} + (k+2\l-1)b_{n-2} + \l b_{n-3}.&
\end{flalign}
\end{gather}
Moreover, the sequence $b_n$ satisfies a shorter recurrence if and only
if $k = 0$, $ k=1 $, or $\l = 0$. When $k = 0$, this recurrence is given by $b_n = \l b_{n-1}$, and 
when $k = 1$, it is
\begin{equation}
b_{n} =  \l b_{n-1} + \l b_{n-2}.
\end{equation}
\end{lemma}

\begin{proof}
The first recurrence follows from deleting an outer edge of the path graph $P_n$
and using Lemma \ref{edgerecursion}.  To verify the second one, we first use
Lemma \ref{edgerecursion} (picking any edge in $C_n$) to give
\begin{equation}\label{bnanrecur}
b_n = a_{n} - b_{n-1} + \l a_{n-2}.
\end{equation}
Let $c_n = b_n + b_{n-1} = a_{n} + \l a_{n-2}$ and notice that $c_n$ satisfies
the same recurrence as $a_{n}$; namely,
\begin{equation}
\begin{split}
c_n =   \ & a_{n}+\l a_{n-2}  \\ 
= \ & (k+\l-1)a_{n-1}+ \l a_{n-2}+\l\left((k+\l-1)a_{n-3}+ \l a_{n-4}\right)  \\
= \ &  (k+\l-1)(a_{n-1}+\l a_{n-3}) + \l(a_{n-2}+\l a_{n-4})  \\  
= \ & (k+\l-1)c_{n-1} + \l c_{n-2}. \\ 
\end{split}
\end{equation}
It follows that $b_n$ satisfies the third order recurrence given in the statement
of the lemma.  Additionally, the initial conditions for both sequences $a_{n}$ and $b_n$
are easily worked out to be the ones shown.

Finally, suppose that the sequence $b_n$ satisfies a shorter recurrence, 
\[b_{n} + r b_{n-1} + s b_{n-2} = 0,\] and let 
\begin{equation*}
B = \left[ \begin{array}{ccc}b_3 & b_2 & b_1 \\b_4 & b_3 & b_2 \\b_5 & b_4 & b_3\end{array}\right].
\end{equation*}
Since the nonzero vector $[1,r,s]^T$ is in the kernel of $B$, we must have that
\begin{equation*}
0 = \det(B) = -k^2(k-1)\l((k+\l-1)^2+4\l).
\end{equation*}
It follows that for $b_n$ to satisfy a smaller recurrence, we must have 
$k = 0$, $ k=1 $, or $\l = 0$.  It is clear that when $k = 0$, we have
$b_n = \l^n = \l b_{n-1}$.  When $k = 1$, we can use Lemma \ref{vertexrecursion} to see that
\[ b_{n+1} = \l (a_{n}+\l a_{n-2}),\]
and combining this with (\ref{bnanrecur}) gives the recurrence stated in the lemma.
\end{proof}

When $k = 1$ and $\l = 1$, the recurrences given by Lemma \ref{recurlem} when applied to the families of path graphs and cycle graphs are the Fibonacci and Lucas numbers, respectively. This observation is well-known (see \cite[Examples 4.1 and 5.3]{Koshy01}) and was brought to our attention by Cox \cite{C}:
\begin{equation}\label{fiblucas}
\chi_{P_n}(1,1) = F_{n+2} \qquad \text{and} \qquad \chi_{C_n}(1,1) = L_{n}.  
\end{equation}

Moreover, when $ k=2 $ and $ \l=1 $, the recurrence given by Lemma \ref{recurlem} when applied to the family of path graphs is the one associated to the Pell numbers:
\[
\chi_{P_n}(2,1) = Q_{n+1},
\]
where $Q_0 = 1$, $Q_1 = 1$, and $Q_n = 2Q_{n-1}+Q_{n-2}$.

\section{Identities}

In this section, we derive some identities involving 
the generalized Fibonacci and Lucas numbers $ a_n $ and $ b_n $ using the graph coloring interpretation found here. In what follows, we fix $ k=1 $. In this case, the $ a_n $ and $ b_n $ 
satisfy the following recurrences:
\[a_n = \l a_{n-1}+\l a_{n-2} \qquad \text{and} \qquad b_n = \l b_{n-1}+\l b_{n-2}.\]

\begin{theorem}\label{main}
The following identities hold:
\begin{align}
b_n & = \l a_{n-1}+\l^2 a_{n-3},\\
b_n & = a_{n}-\l^2 a_{n-4},\\
a_{r+s} & = \l a_ra_{s-1}+\l^2a_{r-1}a_{s-2},\\
a_{r+s} & = a_ra_s-\l^2a_{r-2}a_{s-2},\label{eqn4}\\
a_{r+s+t+1} & = \l a_ra_sa_t + \l^3a_{r-1}a_{s-1}a_{t-1}-\l^4a_{r-2}a_{s-2}a_{t-2}.
\end{align}
\end{theorem}

\begin{proof}
All the identities in the statement of the theorem follow from \mbox{Lemma \ref{vertexrecursion}} when applied to different graphs (with certain choices of vertices and edges). To see the first two equations, consider the cycle graph $ C_n $ and pick any vertex and any edge. To see the next two equations, consider the path graph $ P_{r+s} $ with $ v = r+1 $ and $ e = \set{r,r+1} $.
\begin{center}
\includegraphics{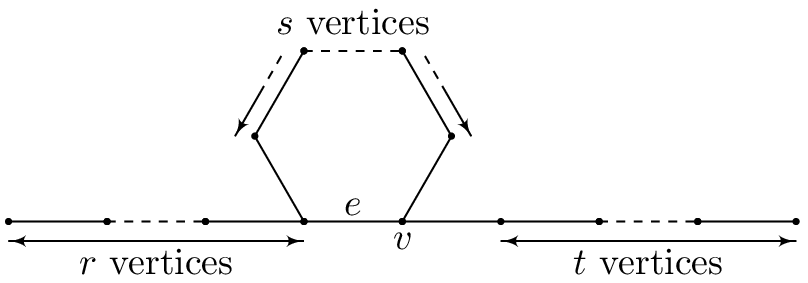}
\end{center}
In order to prove the final equation in the statment of the theorem, consider the graph $ G $ in the above figure. It follows from Lemma \ref{vertexrecursion} that
\[
\l a_{r+s}a_t + \l^3 a_{r-1}a_{s-1}a_{t-1} = a_{r+s+t+1} - \l^4 a_{r-2}a_{s-2}a_{t-1}.
\]
Rearranging the terms and applying \eqref{eqn4}, we see that
\begin{align*}
a_{r+s+t+1} = & \ \l a_{r+s}a_t + \l^3 a_{r-1}a_{s-1}a_{t-1} + \l^4 a_{r-2}a_{s-2}a_{t-1}\\
 = & \ \l (a_ra_s-\l^2a_{r-2}a_{s-2})a_t + \l^3 a_{r-1}a_{s-1}a_{t-1} + \l^4 a_{r-2}a_{s-2}a_{t-1}\\
 = & \ \l a_ra_sa_t -\l^3a_{r-2}a_{s-2}(\l a_{t-1}+\l a_{t-2})\\
& \qquad + \l^3 a_{r-1}a_{s-1}a_{t-1} + \l^4 a_{r-2}a_{s-2}a_{t-1}\\
 = & \ \l a_ra_sa_t + \l^3a_{r-1}a_{s-1}a_{t-1}-\l^4a_{r-2}a_{s-2}a_{t-2}.
\end{align*}
This completes the proof of the theorem.
\end{proof}

\begin{corollary}
The following identities hold:
\begin{align*}
L_n&=F_{n+1}+F_{n-1},\\
L_n&=F_{n+2}-F_{n-2},\\
F_{r+s} &= F_{r+1}F_{s}+F_{r}F_{s-1},\\
F_{r+s} &= F_{r+1}F_{s+1}-F_{r-1}F_{s-1},\\
F_{r+s+t} &= F_{r+1}F_{s+1}F_{t+1}+F_rF_sF_t-F_{r-1}F_{s-1}F_{t-1}.
\end{align*}
\end{corollary}

\begin{proof}
The indentities follow from the corresponding ones in Theorem \ref{main} with $\l=1$ 
by making suitable shifts of the indices and using \eqref{fiblucas}.
\end{proof}

\section{Further Exploration}

In this note, we have produced recurrences and identities 
by decomposing different classes of graphs in different ways.  Our treatment
is by no means exhaustive, and there should be many ways to expand on what
we have done here.  For instance, is there a graph coloring proof of Cassini's identity?

%---------------------------------------------------------------------------
\bibliographystyle{amsalpha}%{alpha}%{plain} %{amsplain} 
\bibliography{fibcolorings}

\end{document}